
\mag=\magstep1
\documentstyle{amsppt}
\input amsppt1
\pageheight{23.5true cm}
\pagewidth{15.5true cm}
\parindent=4mm
\baselineskip=13pt plus.1pt
\parskip=3pt plus1pt minus.5pt
\nologo
\NoRunningHeads
\NoBlackBoxes

\topmatter

\title Flips of semi-stable 4-folds
       whose degenerate fibers are unions of Cartier divisors
       which are terminal factorial 3-folds
\endtitle
\author Yasuyuki Kachi \endauthor
\affil Department of Mathematical Sciences \\
              University of Tokyo \\
            Hongo, Tokyo 113, Japan \\
   E-mail : kachi\@ \! \! clover.math.s.u-tokyo.ac.jp \endaffil

\abstract
{We shall investigate flipping contractions from
a semi-stable 4-fold $X$ whose degenerate fiber is a union
of Cartier divisors which are terminal factorial 3-folds.
Especially we shall prove that $X$ is smooth along the
flipping locus, and that the flip exists for such
contractions.}

\toc
\subhead
 \S 0.\, \, Introduction
 \endsubhead
\subhead
 \S 1.\, \,
  Preliminaries (Dividing into cases)
 \endsubhead
\subhead
 \S 2.\, \,
 The case $E_i \subset D_1 \cap D_2$
 for some $E_i$
 \endsubhead
\subhead
 \S 3.\, \,
 The case \, $E_i \subset D_1$, \,
 $E_i \not\subset \bigcup\limits_{k \geq 2} D_k$ \, and \,
 $E_i \cap \bigcup\limits_{k \geq 2} D_k \not= \emptyset$
 \endsubhead
\subhead \, \,
 for any $E_i$
 \endsubhead
\subhead
 \S 4.\, \,
 The exclusion of the case $E \subset D_1$ and
 $E \cap \bigcup\limits_{k \geq 2} D_k = \emptyset$
 \endsubhead
\subhead
 \S 5.\, \,
 Description of flips
 \endsubhead
\endtoc
\endtopmatter

\document

We will work over $\Bbb C$,
the complex number field.

\head \S 0.\ Introduction. \endhead

The {\it Minimal Model Conjecture\/} asserts that
each algebraic variety has either a minimal model, or a
model which admits a structure of Mori fiber space.
This was classical in dimension 2, while has been
considered to be extremely difficult
in dimension greater than or equal to 3. By virtue
of the works of Reid [R1,3], Mori [Mo2], Kawamata [Kaw1,2]
and Shokurov [Sh1], however, this was reduced to the
{\it existence of flips\/} (see [R2] Conjecture 3.7,
[Kaw2] Problem 5.6).
In the case of semi-stable 3-folds, Kawamata [Kaw3]
solved this affirmatively. His idea is to take
double coverings to reduce the problem to
the existence of {\it flops\/}, much simpler than
that of flips ({\it cf.\/} Koll\'ar [Ko1]).
By applying [Kaw3],
Mori [Mo4] then proved the existence of flips of 3-folds
in the general case, and thus the Minimal Model Conjecture
has been found to be true also in dimension 3.
There are also further developments on 3-dimensional flips,
such as [Utah], [KoMo], [Sh2], and [Kaw5]
(see also [R5]).

Thus it is worth trying to investigate flips in
dimension greater than or equal to 4 as well. In this direction,
little was known, except Kawamata's structure theorem [Kaw4]
in the case of smooth 4-folds:

\proclaim{Theorem 0.1}\ \ (Kawamata [Kaw4])

\flushpar
Let $X$ be a smooth projective 4-fold and
$g : X \to Y$ a flipping contraction. Then
$\text{Exc }g$ is a disjoint union of $\Bbb P^2\text{'s}$,
and those normal bundles in $X$ are all
isomorphic to $\Cal O_{\Bbb P^2}(-1)^{\oplus 2}$.

Moreover, the flip of $g$ exists.
\endproclaim

This, together with the {\it Termination Theorem\/} [KaMaMa] after
[Sh2], should be considered as the first step of generalizing the
Minimal Model Conjecture to dimension 4.

In this paper, we shall investigate the structures of
flipping contractions from semi-stable 4-folds,
with a certain additional assumption
(Main Theorem 0.5).
This, together with Kawamata's Theorem [Kaw4]
(see also Theorem 0.1 above), shows that
the flip exists for such contractions.

\definition{Assumption 0.2}\ \ (Semi-stable degenerations)
\, ({\it cf.\/} [Kaw3,5])

Let $f : X \to \Delta$ be a projective morphism from a 4-dimensional
analytic space $X$ with at most terminal singularities to the
disc $\Delta := \{ z \in \Bbb C \, | \, |z|<1\}$ such that
each fiber of $f$ over $\Delta - \{0\}$ is a projective
3-fold with at most terminal singularities.
Moreover, assume the following conditions:

\flushpar
(0.2.0) \, \, \, The central fiber $D := f^{-1}(0)$ is reduced,

\flushpar
(0.2.1) \, \, \, $D$ is normal crossing outside $\text{Sing } X$,

\flushpar
(0.2.2) \, \, \,
Let $D = \bigcup\limits_{k=1}^r D_k$ be the
irreducible decomposition, then each $D_k$ is a normal 3-fold which
is a $\Bbb Q\text{-Cartier}$ divisor on $X$, and

\flushpar
(0.2.3) \, \, \,
The pair $(X \, , \, D)$ is log-terminal.
\enddefinition

\definition{Definition 0.3}\ \
Let $f : X \to \Delta$ be as in Assumption 0.1.
Then the contraction morphism
$g : X \to Y$ associated to an extremal
ray $R$ of $\overline{NE}(X/\Delta)$
is said to be a {\it flipping contraction\/}
if $\dim \text{Exc }g \leq 2$, and
$\text{Exc } g$ is called the {\it flipping locus\/}.
If there is a projective bimeromorphic morphism $g^+ : X^+ \to Y$
over $\Delta$ from another 4-dimensional analytic space $X^+$
such that

\flushpar
(0.3.1) \, \, \,
The composite $X^+ \overset {g^+}\to\longrightarrow Y
\longrightarrow \Delta$
satisfies the Assumption 0.2,
except possibly (0.2.1),

\flushpar
(0.3.2) \, \, \,
$\rho(X^+/Y)=1$,

\flushpar
(0.3.3) \, \, \,
$K_{X^+}$ is $g^+\text{-ample}$, and

\flushpar
(0.3.4) \, \, \,
$\dim \text{Exc }g^+ \leq 2$,

\flushpar
then $g^+$ is called the {\it flip\/}
of $g$. We sometimes call the composite bimeromorphic map
$g^{+ \, -1} \circ g : X \dashrightarrow X^+$
also the {\it flip\/} of $g$.
(As for the condition (0.3.1), see Remark 5.4.)
\enddefinition

\definition{Remark 0.4}\ \
Note that the flip $g^+$ of $g$ is unique if exists.
\enddefinition

\proclaim{Main Theorem 0.5}\ \
In addition to the Assumption 0.2 and Definition 0.3,
we assume furthermore the followings:

\flushpar
(0.5.1) \, \, \, Each $D_i$ is a Cartier divisor
which is a terminal factorial 3-fold.

Then, $\text{Exc } g$ is a disjoint union of
$\Bbb P^2\text{'s}$, and $X$ is smooth along
$\text{Exc } g$. Let us fix a connected
component $E \simeq \Bbb P^2$ of $\text{Exc } g$
arbitrarily. Then
$$
N_{E/X} \simeq \Cal O_{\Bbb P^2}(-1)^{\oplus 2}
$$
and after a suitable renumbering of $\{ D_1 , \dots , D_r\}$,
exactly one of the followings holds:

(A-I) \, \, \,
$E$ is a connected component of $D_1 \cap D_2$,
$N_{E/D_k} \simeq \Cal O_{\Bbb P^2}(-1)$
$(k = 1,2)$, $D_3 \cap E$, $D_4 \cap E$ are both lines
in $E$ which are distinct to each other, and
$D_k \cap E = \emptyset$ $(\forall \, k \geq 5)$,

(A-II) \, \, \,
$E$ is a connected component of $D_1 \cap D_2$,
$N_{E/D_k} \simeq \Cal O_{\Bbb P^2}(-1)$
$(k = 1,2)$, $D_3 \cap E$ is a smooth conic
in $E$, and
$D_k \cap E = \emptyset$ $(\forall \, k \geq 4)$.

or

(B) \, \, \,
$E \subset D_1$, $N_{E/D_1} \simeq \Cal O_{\Bbb P^2}(-1)$,
$D_2 \cap E$ is a line in $E$, and
$D_k \cap E = \emptyset$ $(\forall \, k \geq 3)$.

Moreover, the flip $g^+$ of $g$ exists.
\endproclaim

As for the description of the flip $g^+$, see \S 5.

\proclaim{Remark 0.6}\ \

\flushpar
(1) \, \, \,
There are no flipping contractions from
terminal Gorenstein 3-folds [Mo4].

\flushpar
(2) \, \, \,
In particular, we may assume that $E \subset D$,
in Main Theorem 0.5.
\endproclaim

\definition{Notation 0.7}\ \
Let $R$ be the extremal ray of $\overline{NE} (X/\Delta)$
determining the flipping contraction $g : X \to Y$.
Then for $D' = \sum\limits_{k=1}^r \alpha_k D_k$
we simply write $(D' \, . \, R) >0$
({\it resp.\/} $\geq 0$, $=0$, $\leq 0$, $<0$)
when for an irreducible curve $C \subset E$,
$(D' \, . \, C)>0$ ({\it resp.\/}
$\geq 0$, $=0$, $\leq 0$, $<0$).
In particular for $D = \sum\limits_{k = 1}^r D_k$,

\flushpar
(0.7.1) \qquad \qquad \qquad \qquad \qquad \qquad \, \, \,
$(D \, . \, R)=0$.
\enddefinition

\definition{Acknowledgement}\ \
The author would like to express his gratitude for
Professor Y.Kawamata who suggested to him a
simplified proof of 2.2 in the case that at least one
irreducible component of the flipping surface is
contained in the double locus of the degenerate fiber.

He is grateful to Professor J. Koll\'ar who gave him
lots of helpful advices during his stay at University
of Utah in April, 1995, especially suggested to him
a simplified proof 4.3 for excluding the case that the flipping
surface is disjoint from the double locus of the degenerate fiber.

He is grateful to Professor S.Mukai who told him
a related example, which inspired him to this problem.

He is grateful to Professor S.Mori
who gave him some linguistic advices.

He also would like to thank Professors S.Ishii and
K.-i.Watanabe for their helpful advices, and
Professors D.Burns, A.Bertram, Y.Ruan, T.Maeda and K.Oguiso
for their stimulating conversations and encouragements.
\enddefinition

\head \S 1.\ Preliminaries (Dividing into cases). \endhead

\proclaim{Lemma 1.1}\ \
$E$ is purely 2-dimensional.
\endproclaim

\demo{Proof}\ \
Assume that $E$ has an irreducible component $C$
of dimension 1 to derive a contradiction. First we claim:

\flushpar
(1.1.1) \, \, \,
There is at least one $D_k$ which contains $C$ such that
$(-K_{D_k} \, . \, C) >0$.

Actually, assume $D_1$ say, satisfies
$D_1 \supset C$ and $(-K_{D_1} \, . \, C)\leq 0$.
Since $(D \, . \, C)=0$,
$$
\align
(1.1.2) \qquad \qquad \quad \,
0 \geq (-K_{D_1} \, . \, C)_{D_1}
&= (-K_X \, . \, C)_X + (-D_1 \, . \, C)_X
\qquad \qquad \qquad \quad \, \, \, \, \,
\\
&\geq 1 + (-D_1 \, . \, C)_X \\
&= 1 + \biggl( \sum\limits_{k \geq 2} D_k \, . \, C \biggr)_X \, .
\endalign
$$
In particular

\flushpar
(1.1.3) \, \, \, $(D_k \, . \, C)_X <0$ for some $i \geq 2$.

\flushpar
For such $k$, $C \subset D_k$ and in the same way as in (1.1.2),
we have
$$
\align
(1.1.4) \qquad \qquad \qquad \,
(-K_{D_k} \, . \, C)_{D_k}
&= (-K_X \, . \, C)_X + (-D_k \, . \, C)_X
\qquad \qquad \qquad \qquad \, \, \, \,
\\
& >0 \qquad (\text{by } (1.1.3))
\endalign
$$
and hence (1.1.1).

Let $V$ be an analytic neighborhood of $g(D_k \cap E)$
in $g(D_k)$, $U := g^{-1}(V)$, and
$\widetilde{V}$ the normalization of $V$.
Let $(g|_{U})\sptilde : U \to \widetilde{V}$
be the morphism induced from
$g|_{U} : U \to V$.
(Recall that $D_k$ is assumed to be normal (0.2.2).)
This is a bimeromorphic morphism
with $\text{Exc } (g|_U)\sptilde = D_k \cap E$, and
$C$ forms an irreducible component of
$\text{Exc }(g|_U)\sptilde$.
Let $C^-$ be the union of all the other irreducible components
of $\text{Exc }(g|_U)\sptilde$.
Let $L + L^-$ be a
$(g|_U)\sptilde\text{-very}$ ample divisor of $U$
such that $L \cap C =\emptyset$ and $L^- \cap C^- =\emptyset$.
Then $|L|$ determines
the morphism $h : U \to U'$ over $\widetilde{V}$
which is also a bimeromorphic morphism such that
$$
\text{Exc }h = C.
$$
Since $(-K_{D_k} \, . \, C)>0$ (1.1.1),
this is a flipping contraction from $U$.
Since $D_k$ is assumed to be a terminal Gorenstein 3-fold
(0.1.3), this contradicts Remark 0.6.
Hence the Lemma. \quad \qed
\enddemo

\proclaim{Lemma 1.2}\ \
Let $C$ be a smooth rational curve which is
contained in $E$. Assume that
$C \cap \text{Sing }X = \emptyset$.
Let $\text{Hilb}_{X, [C]}$ be
the connected component of the Hilbert scheme
$\text{Hilb}_X$ containing the point $[C]$.
Then
$$
\dim \text{Hilb}_{X, [C]} \geq 2.
$$
\endproclaim

\demo{Proof}\ \
Since $\dim X =4$ and $(-K_X \, . \, C) >0$,
this is a direct consequence of [Mo1]:
$$
\dim \text{Hilb}_{X, [C]} \geq \text{codim}_X C
+ (-K_X \, . \, C) + \text{deg } K_C \geq 2.
\quad \qed
$$
\enddemo

\proclaim{Corollary 1.3}\ \
Each $E_i$ has no open subset which admits a
$\Bbb P^1\text{-bundle}$ structure.

In particular, $g(E)$ is a point. \quad \qed
\endproclaim

\proclaim{Proposition 1.4}\ \ (Dividing into cases)

After a suitable renumbering of $\{D_1, \dots , D_r\}$,
exactly one of the followings holds:

\flushpar
(A) \, \, \, $E_i \subset D_1 \cap D_2$, \,
and \, $E_i \not\subset \bigcup\limits_{k \geq 3}D_k$
\, for some $E_i$,

\flushpar
(B) \, \, \, $E_i \subset D_1$, \,
$E_i \not\subset \bigcup\limits_{k \geq 2} D_k$,
\, and \, $E_i \cap \bigcup\limits_{k \geq 2} D_k
\not= \emptyset$ \, for all $E_i$, \, or

\flushpar
(C) \, \, \, $E \subset D_1$, \,
and \, $E \cap \bigcup\limits_{k \geq 2} D_k = \emptyset$.
\endproclaim

\demo{Proof}\ \
Assume not the Case (A). Then by the assumption (0.2.1)
with Lemma 1.1, for any $E_i$
there exists some $D_{k(i)}$ such that
$$
E_i \subset D_{k(i)}, \, \text{ and } \,
E_i \not\subset \bigcup\limits_{k \not= k(i)} D_k.
$$
In particular
$\Bigl( \sum\limits_{k \not= k(i)} D_k \, . \, R \Bigr)
\geq 0$, or equivalently
$$
(D_{k(i)} \, . \, R) \leq 0.
$$
If $(D_{k(i)} \, . \, R) <0$,
then
$$
E \subset D_{k(i)} \text{ and } \,
E \cap \bigcup\limits_{k \not= k(i)}
D_{k(i)} \not= \emptyset,
$$
thus we have the Case (B) by letting
$D_1 = D_{k(i)}$.
Otherwise we have the
Case (C).
\quad \qed
\enddemo

The following is due to Cutkosky [C]
which is a generalization of Mori [Mo2], and
is an application of Fujita [F]:

\proclaim{Theorem 1.5}\ \ (Cutkosky [C])

Let $Z$ be a terminal factorial 3-fold and
$h: Z \to W$ the contraction of an extremal ray
of $Z$ which is birational. Let
$F := \text{Exc }h$ and assume that
$h(F)$ is a point. Then
$(F, \Cal O_F(F))$ is isomorphic either to
$$
(\Bbb P^2, \Cal O_{\Bbb P^2}(-1)),
\, \,
(\Bbb P^2, \Cal O_{\Bbb P^2}(-2)),
\, \,
(\Bbb P^1 \times \Bbb P^1, \Cal O_{\Bbb P^1 \times \Bbb P^1}(-1,-1)),
\, \,
\text{ or }
\, \, (S_2, \Cal O_{S_2}(-1)),
$$
where $S_2$ is the singular quadrtic surface in $\Bbb P^3$.
$h$ is the blow-up with the reduced center $P := h(F)$.

Furthermore, in the case
$(F, \Cal O_F(F)) \simeq (\Bbb P^2, \Cal O_{\Bbb P^2}(-1))$
$(${\it resp.\/} $(\Bbb P^2, \Cal O_{\Bbb P^2}(-2))$, we have

\flushpar
(1) \, \, \, $Z$ is smooth along $F$, and

\flushpar
(2) \, \, \,
$W$ is smooth at $P$
$(${\it resp.\/} $W$ has a quotient singularity
of type $\dfrac{1}{2}(1,1,1)$ at $P)$.
\endproclaim

\head \S 2.\ The case $E_i \subset D_1 \cap D_2$ for
some $E_i$. \endhead

The aim of this section is to prove the following:

\proclaim{Theorem 2.1}\ \
In the Case (A) in Proposition 1.4, the followings hold:

\flushpar
(1) \, \, \, $E$ is irreducible: $E = E_i$, and is
isomorphic to $\Bbb P^2$.

\flushpar
(2) \, \, \, $X$ is smooth along $E$.

\flushpar
(3) \, \, \, $N_{E/D_1} \simeq N_{E/D_2} \simeq
\Cal O_{\Bbb P^2}(-1)$. In particular
\, $N_{E/X} \simeq \Cal O_{\Bbb P^2} (-1)^{\oplus 2}$.

\flushpar
(4) \, \, \, After a suitable renumbering of $\{D_3, \dots , D_r\}$,
either one of the followings holds:

(A-I) \, \, \,
$D_3 \cap E$, $D_4 \cap E$ are distinct lines in $E$, and
$D_k \cap E = \emptyset$ $(\forall k \geq 5)$,

or

(A-II) \, \, \,
$D_3 \cap E$ is a smooth conic in $E$, and
$D_k \cap E = \emptyset$ $(\forall k \geq 4)$,

\flushpar
(5) \, \, \, The flip $g^+$ of $g$ exists.
\endproclaim

First we shall prove the following, which is told us by
Y.Kawamata and J.Koll\'ar:

\proclaim{Proposition 2.2}\ \ (Kawamata - Koll\'ar)

In the Case (A), $E$ is irreducible, and hence is
a connected component of $D_1 \cap D_2$.
\endproclaim

\demo{Proof}\ \
Assume that $E$ is reducible, and let $E_1$
be an irreducible component
of $E$ such that

\flushpar
(2.2.1) \qquad \qquad \qquad \qquad \qquad \qquad
$E_1 \subset D_1 \cap D_2$.

\flushpar
Since $E$ is connected,
there exists another irreducible component, say $E_2$,
of $E$ such that

\flushpar
(2.2.2) \qquad \qquad \qquad \qquad
$E_2 \not\subset D_1 \cap D_2$ \, and \,
$E_1 \cap E_2 \not= \emptyset$.

\flushpar
In particular, $(D_k \, . \, R)>0$ for either
$k = 1 \text{ or } 2$. Let us assume

\flushpar
(2.2.3) \qquad \qquad \qquad \qquad \qquad \qquad \, \,
$(D_2 \, . \, R) >0$.

On the other hand, by (2.2.1) with the assumption (0.2.1)
and Lemma 1.1,
$E_1 \not\subset D_k$ for all $k \geq 3$, thus

\flushpar
(2.2.4) \qquad \qquad \qquad \qquad \qquad \, \,
$(D_k \, . \, R) \geq 0 \, \, \, (\forall k \geq 3)$.

\flushpar
{}From (2.2.3) and (2.2.4), we necessarily have

\flushpar
(2.2.5) \qquad \qquad \qquad \qquad \qquad \qquad \, \,
$(D_1 \, . \, R) <0$,

\flushpar
and in particular

\flushpar
(2.2.6) \qquad \qquad \qquad \qquad \qquad \qquad \quad \, \,
$E \subset D_1$.

\flushpar
(2.2.7) \, \, \,
Let $g|_{D_1} : D_1 \to g(D_1)$ be the restriction of $g$
to $D_1$, let $g(D_1)\sptilde$ be the normalization of
$g(D_1)$, and $h : D_1 \to g(D_1)\sptilde$
the morphism induced from $g|_{D_1}$.

Then $h$ is a birational morphism such that
$$
\text{Exc }h = E.
$$
Thus for a general irreducible curve $C$ in $E_1$,

\flushpar
(2.2.8) \qquad \qquad \qquad \qquad \qquad \quad \, \, \, \, \,
$(E_1 \, . \, C)_{D_1} <0$.

\flushpar
This contradicts (2.2.3), since
$(E_1 \, . \, C)_{D_1} = (D_2 \, . \, C)_X$ (2.2.1).
Hence $E$ must be irreducible.
\quad \qed
\enddemo

\proclaim{Proposition 2.3} \qquad \qquad \qquad \qquad
$E \simeq \Bbb P^2$.
\endproclaim

\demo{Proof}\ \

\flushpar
(2.3.0) \, \, \,
Let $g(D_k)\sptilde$ be the normalization of
$g(D_k)$ $(k = 1,2)$, and let
$h_k : D_k \to g(D_k)\sptilde$
be the morphism induced from $g|_{D_k}$.

First we claim that

\flushpar
(2.3.1) \qquad \quad
$(D_k \, . \, R) <0$ \, and \, $-K_{D_k}$ is $h_k\text{-ample}$
\, for both \, $k=1,2$.

In fact, let $C$ be a general irreducible curve in $E$.
Then as in (2.2.8)
$$
(D_2 \, . \, C)_X = (E \, . \, C)_{D_1} <0,
$$
and similarly $(D_1 \, . \, C)_X <0$.
Thus
$$
(-K_{D_k} \, . \, C)_{D_k}
= (-K_X \, . \, C)_X - (D_k \, . \, C)_X
>0 \, \, (k = 1,2),
$$
and we get (2.3.1).

By $\text{Exc } h_k = E$, Corollary 1.3 and Proposition 2.2,
$$
\rho(D_k/g(D_k)\sptilde) =1
$$
and it follows from [C], together with
Corollary 1.3, that

\flushpar
(2.3.2) \qquad \qquad \qquad \qquad \qquad \qquad
$E \simeq \Bbb P^2$ \, or \, $S_2$.

\flushpar
Assume $E \simeq S_2$ to get a contradiction.

Let $B$ be a general $(+2)\text{-section}$ of
$E \simeq S_2$. Since $D_1$ is a Cartier divisor
of $X$ and since $E$ is a Cartier divisor of $D_1$
(0.5.1), $X$ is smooth along $E$ outside the vertex.
Hence

\flushpar
(2.3.3) \qquad \qquad \qquad \qquad \qquad \quad \, \, \, \,
$B \cap \text{Sing } X = \emptyset$.

\flushpar
Consider the exact sequence:
$$
0 \longrightarrow N_{B/E} \longrightarrow
N_{B/X} \longrightarrow N_{E/X}
\otimes \Cal O_B \longrightarrow 0
$$
Since $E$ is a connected component of $D_1 \cap D_2$,
$$
N_{E/X} \otimes \Cal O_B \simeq
\Cal O_B (D_1) \oplus \Cal O_B (D_2),
$$
and we have
$$
\align
(2.3.4) \qquad \qquad \qquad
c_1(N_{B/X}) &= c_1(N_{B/E})
+ (D_1 \, . \, B) + (D_2 \, . \, B)
\qquad \qquad \qquad \qquad \, \, \, \, \\
&= 2 + (D_1 \, . \, B) + (D_2 \, . \, B).
\endalign
$$

On the other hand, since
$B \cap \text{Sing } X = \emptyset$ (2.3.3),
we can consider the deformation of $B$ inside $X$
[Mo1] ({\it cf.\/} [I], [Wi\'s], [Ko2]).
Let $T:= \text{Hilb}_{X, \, [B]}$ be the
connected component of the Hilbert scheme
$\text{Hilb}_X$ containing the point
corresponding to $B$. Then
$$
\align
(2.3.5) \qquad \qquad \qquad \quad \, \, \,
3 = \dim T &\geq \dim X + (-K_X \, . \, B) - 3
\qquad \qquad \qquad \qquad \qquad \, \, \\
&= (-K_X \, . \, B) + 1.
\endalign
$$

Moreover, since $-K_X$ is Cartier (0.5.1), and since
the numerical class of $B$ is the double of the class
of a ruling of $E$,
$$
(-K_X \, . \, B) \geq 2.
$$
Hence the inequality (2.3.5) must be the equality:

\flushpar
(2.3.6) \qquad \qquad \qquad \qquad \qquad \quad \, \, \, \, \,
$(-K_X \, . \, B) = 2$,

\flushpar
that is,

\flushpar
(2.3.7) \qquad \qquad \qquad \qquad \qquad \qquad
$c_1(N_{B/X})=0$.

{}From this and (2.3.4),
$$
(D_1 \, . \, B) + (D_2 \, . \, B) = -2.
$$
Again since $(D_k \, . \, B)$ is an even integer
$(k = 1,2)$,

\flushpar
(2.3.8) \qquad \qquad \qquad \quad \, \,
$(D_k \, . \, B) \geq 0 \, \,
\text{ for either } \, k = 1 \text{ or } 2$,

\flushpar
which contradicts (2.3.1).
Hence we must have
$E \simeq \Bbb P^2$ (2.3.2).
\quad \qed
\enddemo

\proclaim{Corollary 2.4}\ \
$X$ is smooth along $E$, and
$$
N_{E/D_1} \simeq N_{E/D_2} \simeq
\Cal O_{\Bbb P^2}(-1).
$$
\endproclaim

\demo{Proof}\ \
First by Proposition 2.3 and Lemma 1.5,
$D_1$, $D_2$ are smooth along
$E$, and hence so is $X$,
since $D_1$ is a Cartier divisor of $X$.

Next, let $h_k : D_k \to g(D_k)\sptilde$ be as in
(2.3.0). Since $\text{Exc } h_k = E$
for both $k=1,2$,
$$
N_{E/D_k} \simeq \Cal O_{\Bbb P^2} (a_k) \, \,
\text{ with } \, a_k <0 \, \, \, (k=1,2).
$$
In particular
$N_{E/X} \simeq \Cal O_{\Bbb P^2} (a_1)
\oplus \Cal O_{\Bbb P^2} (a_2)$.
Moreover since $-K_X$ is $g\text{-ample}$,
we must have

\flushpar
\qquad \qquad \qquad \qquad \qquad \qquad \qquad \quad \,
$a_1 = a_2 = -1$. \quad \qed
\enddemo

\proclaim{Corollary 2.5}\ \
The flip of $g$ exists.
\endproclaim

\demo{Proof}\ \
Since $X$ is smooth along $E$,
this is a direct consequence of Kawamata [Kaw4].
\quad \qed
\enddemo

The rest thing we have to prove is the following:

\proclaim{Proposition 2.6}\ \
After a suitable renumbering of $\{D_3, \dots , D_r\}$,
either one of the followings holds:

(I) \, \, \,
$D_3 \cap E$, $D_4 \cap E$ are distinct lines in $E$, and
$D_k \cap E = \emptyset$ $(\forall k \geq 5)$,

or

(II) \, \, \,
$D_3 \cap E$ is a smooth conic in $E$, and
$D_k \cap E = \emptyset$ $(\forall k \geq 4)$.
\endproclaim

\demo{Proof}\ \
Let $l$ be any line in $E$. Since
$N_{E/D_k} \simeq \Cal O_{\Bbb P^2} (-1)$ \, $(k=1,2)$,
$$
\align
(D_k \, . \, l)_X &= (-K_X \, . \, l)_X
+ (K_{D_k} \, . \, l)_{D_k} \\
&= 1-2 = -1 \, \, \,
(k=1,2).
\endalign
$$
Thus either

\flushpar
(2.6.1) \qquad \qquad \qquad \, \, \, \,
$(D_3 \, . \, l) = 2$, \, \, \, $(D_k \, . \, l)=0$
\, $(\forall \, k \geq 4)$, \, \, or

\flushpar
(2.6.2) \qquad \qquad \quad \, \,
$(D_3 \, . \, l) = (D_4 \, . \, l) = 1$, \, \, \, $(D_k \, . \, l)=0$
\, $(\forall \, k \geq 5)$.

In the case (2.6.1), $D_3 \cap E$ must be a smooth conic $C$
in $E \simeq \Bbb P^2$, by the assumption (0.2.1).
Hence we have the case (II).

On the other hand, in the case (2.6.2),
obviously we have (I).
\quad \qed
\enddemo

Now the proof of Theorem 2.1 is completed.

\head \S 3.\ The case $E_i \subset D_1$,
$E_i \not\subset \bigcup\limits_{k \geq 2} D_k$ and
$E_i \cap \bigcup\limits_{k \geq 2} D_k \not= \emptyset$
for any $E_i$. \endhead

The aim of this section is:

\proclaim{Theorem 3.1}\ \
In the Case (B) in Proposition 1.4, the followings hold:

\flushpar
(1) \, \, \,
$E$ is irreducible and is isomorphic to $\Bbb P^2$,

\flushpar
(2) \, \, \,
$X$ is smooth along $E$,

\flushpar
(3) \qquad \qquad \qquad \, \, \, \, \,
$N_{E/D_1} \simeq \Cal O_{\Bbb P^2}(-1)$ \, and \,
$N_{E/X} \simeq \Cal O_{\Bbb P^2}(-1)^{\oplus 2}$,

\flushpar
(4) \, \, \,
After a suitable renumbering of $\{D_2 , \dots , D_r\}$,
$D_2 \cap E$ is a line in $E$, and
\linebreak
$D_k \cap E = \emptyset$
$(\forall \, k \geq 3)$, and

\flushpar
(5) \, \, \,
The flip $g^+$ of $g$ exists.
\endproclaim

First we shall prove:

\proclaim{Lemma 3.2}\ \
In the Case (B),

\flushpar
(1) \qquad \qquad \qquad \qquad \qquad \, \, \, \,
$E \subset D_1$ \, and $(D_1 \, . \, R) <0$.

\flushpar
(2) \, \, \, For any irreducible curve $C$ in $E$,
$$
(-K_{D_1} \, . \, C) \geq 2.
$$
\endproclaim

\demo{Proof}\ \
{\it (1)\/} \, \, \, By the condition (B),
$E \subset D_1$ and
$\biggl( \sum\limits_{k \geq 2} D_k \, . \, R \biggr) >0$,
that is,
$$
(D_1 \, . \, R) <0.
$$

\flushpar
{\it (2)\/} \, \, \,
Let $C \subset E$ be any irreducible curve. Then
$$
(-K_{D_1} \, . \, C) = -(D_1 \, . \, C) + (-K_X \, . \, C)
\geq 2,
$$
since $D_1$ and $-K_X$ are both Cartier.
\quad \qed
\enddemo

\proclaim{Lemma 3.3}\ \
There exists an irreducible component $E_1$, say, of $E$
which is isomorphic to $\Bbb P^2$ such that
$$
N_{E_1/D_1} \simeq \Cal O_{\Bbb P^2}(-1), \, \,
N_{E_1/X} \simeq \Cal O_{\Bbb P^2}(-1)^{\oplus 2},
\, \, \text{ and } \, \,
\Cal O_{E_1}(D_1) \simeq \Cal O_{\Bbb P^2}(-1).
$$
\endproclaim

\demo{Proof}\ \

\flushpar
(3.3.0) \, \, \,
Let $g(D_1)\sptilde$ be the normalization of $g(D_1)$,
and consider the morphism $h : D_1 \to g(D_1)\sptilde$
induced from $g|_{D_1}$, as in (2.2.7), (2.3.0).

By Lemma 3.2 (2), $-K_{D_1}$ is $h\text{-ample}$.
Take the contraction morphism $D_1 \to V$ associated to
any extremal ray of $\overline{NE}(D_1/g(D_1)\sptilde)$.
Then by Cutkosky [C], together with Corollary 1.3 and
Lemma 3.2 (2), the exceptional locus, say $E_1$, must be
isomorphic to $\Bbb P^2$, and

\flushpar
(3.3.1) \qquad \qquad \qquad \qquad \qquad \quad
$N_{E_1/D_1} \simeq \Cal O_{\Bbb P^2}(-1)$.

\flushpar
In particular $D_1$ and hence $X$ is smooth along $E_1$.

Let $\Cal O_{E_1}(D_1) \simeq \Cal O_{\Bbb P^2}(a)$
$(a \in \Bbb Z)$. Then from (3.3.1) and the exact sequence:
$$
0 \longrightarrow N_{E_1/D_1} \longrightarrow N_{E_1/X}
\longrightarrow \Cal O_{E_1}(D_1) \longrightarrow 0
$$
we have
$$
N_{E_1/X} \simeq \Cal O_{\Bbb P^2}(-1) \oplus \Cal O_{\Bbb P^2}(a).
$$
Since $E_1$ never deforms inside $X$, and since
$\Cal O_{E_1}(-K_X)$ is ample, we necessarily have
$a = -1$, and we are done.
\quad \qed
\enddemo

\definition{3.4}\ \
Consider the local flip (Kawamata flip) [Kaw4]
$$
\eta^{(0)} : X \dashrightarrow X^{(1)}
$$
with respect to $E_1$, and let
$$
g^{(1) \, '} : X^{(1)} \to Y
$$
be the structure morphism.
Let $E_1^{+(1)} \simeq \Bbb P^1$ be the flipped curve, and let
$E_i^{(1)}$ $(i \geq 2)$, $D_k^{(1)}$ be the proper transform
of $E_i$, $D_k$ in $X^{(1)}$, respectively.
Then

\flushpar
(3.4.1) \qquad \qquad \qquad \qquad \quad \,
$X^{(1)}$ is smooth along $E_1^{+(1)}$

\flushpar
[loc.cit], and

\flushpar
(3.4.2) \qquad \qquad \, \, \, \,
$E_1^{+(1)} \subset D_2^{(1)}$, \, and \,
$\text{Exc } g^{(1) \, '} = E_1^{+(1)}
\cup \bigcup\limits_{i \geq 2} E_i^{(1)}$.

\flushpar
Hence $g^{(1) \, '}$ is factored through
$$
g^{(1)} : X^{(1)} \to Y^{(1)}
$$
(see the diagram (3.5.1) below) such that

\flushpar
(3.4.3) \qquad \quad \, \, \,
$-K_{X^{(1)}}$ is $g^{(1)}\text{-ample}$ and
$E^{(1)} := \text{Exc } g^{(1)}
= \bigcup\limits_{i \geq 2} E_i^{(1)}$.
\enddefinition

\definition{3.5} \ \
If $\text{Exc }g^{(1)}$ is still reducible, then
return back to the situation of Lemma 3.3, with the
substitutions of $X$, $D_1$ by $X^{(1)}$, $D_1^{(1)}$,
respectively. Then we again find an $E_2^{(1)} \simeq \Bbb P^2$,
say, which has the normal bundle
$\Cal O_{\Bbb P^2}(-1)^{\oplus 2}$ in $X^{(1)}$.
Do the same procedure as in 3.4 above
for $X^{(1)}$ instead of $X$, to get
$X^{(2)}$ and $g^{(2)} : X^{(2)} \to Y^{(2)}$,
satisfying the similar condition to (3.4.3).
If we continue the processes successively
until $\text{Exc }g^{(p)}$ becomes irreducible,
we get the diagram
$$
\CD
@. @. X @. \quad \overset \eta^{(0)}
\to{\dashrightarrow} \quad @.
X^{(1)} @. \quad \overset \eta^{(1)}
\to{\dashrightarrow} \quad @.
X^{(2)} @. \quad \overset \eta^{(2)}
\to{\dashrightarrow} \quad @.
\dots @. \quad \overset \eta^{(p-1)}
\to{\dashrightarrow} \quad @.
X^{(p)} \\
@. @. @V{g}VV @. @V{g^{(1)}}VV @. @V{g^{(2)}}VV
@. @. @. @VV{g^{(p)}}V \\
\Delta @. \quad \longleftarrow \quad
@. Y @. \longleftarrow @.
Y^{(1)} @. \longleftarrow @.
Y^{(2)} @. \longleftarrow @.
\dots @. \longleftarrow @. Y^{(p)}
\endCD
\tag 3.5.1
$$
with

\flushpar
(3.5.2) \quad \, \,
$\text{Exc } g^{(p)} = E_{p+1}^{(p)} \simeq \Bbb P^2$, \,
along which $X^{(p)}$ and $D_1^{(p)}$ are smooth.
\enddefinition

\proclaim{Proposition 3.6} \ \
In Lemma 3.3, $E$ is irreducible: $E = E_1 \simeq \Bbb P^2$.
The flip of $g$ exists.
\endproclaim

\demo{Proof} \ \
First $X^{(p)}$ is smooth
along $E_{p+1}^{(p)}$ (3.5.2). Hence
$X^{(p-1)}$ is smooth along
$$
\eta^{(p-1) \, -1}(E_{p+1}^{(p)} - E_{p+1}^{(p)} \cap E_p^{+ (p)})
= E_{p+1}^{(p-1)} - E_{p+1}^{(p-1)} \cap E_p^{(p-1)},
$$
over which the flip $\eta^{(p-1)}$ is an isomorphism.

On the other hand, $X^{(p-1)}$ is smooth also along
the flipping surface $E_p^{(p-1)} \simeq \Bbb P^2$.
Thus

\flushpar
(3.6.1) \, \, \,
$X^{(p-1)}$ is smooth along $E_p^{(p-1)} \cup E_{p+1}^{(p-1)}
= \text{Exc } g^{(p-1)}$.

If we go further this argument from $X^{(p-1)}$ upstream to
$X^{(0)} = X$, we conclude that

\flushpar
(3.6.2) \qquad \qquad \qquad \qquad \qquad \, \,
$X$ is smooth along $E$.

\flushpar
Hence by Kawamata [Kaw4], $E$ must be irreducible,
$E \simeq \Bbb P^2$, and the flip $g^+$ of $g$ exists.
\quad \qed
\enddemo

\proclaim{3.7}\ \ Proof of Theorem 3.1.
\endproclaim

Since $\Cal O_E(D_1) \simeq \Cal O_{\Bbb P^2}(-1)$
(Lemma 3.3), for any line $l$ in $E$
we have $(D_1 \, . \, l) = -1$, {\it i.e.\/}
$$
\biggl( \sum\limits_{k \geq 2} D_k \, . \, l \biggr) = 1.
$$
Thus we find an unique $D_2$, say, among $\{D_2 , \dots , D_r\}$
such that $D_2 \cap E$ is a line, and $D_k \cap E = \emptyset$
$(\forall \, k \geq 3)$.

Hence we get Theorem 3.1.
\quad \qed

\head \S 4.\ The exclusion of the case $E \subset D_1$ and
$E \cap \bigcup\limits_{k \geq 2} D_k = \emptyset$. \endhead

In this section we shall exclude the Case (C) in Proposition 1.4:

\proclaim{Theorem 4.1}\ \
The Case (C) in Proposition 1.4 never happens.
\endproclaim

We adopt the simplified proof of Theorem 4.1
based on J.Koll\'ar's idea.

\proclaim{Theorem 4.2}\ \ (Koll\'ar - Mori [KoMo] 11.4)

Let $f : U \to \Delta$ be
a surjective morphism from a 4-fold $U$ with at most
terminal singularities to the disc $\Delta$.
Let $U_0$ be its central fiber. Assume that $U_0$ has
at most terminal singularities,
and there is a proper bimeromorphic morphism
$\varphi_0 : U_0 \to V_0$ to a germ $(V_0,Q)$ of
some normal 3-fold such that
$$
R^1 \varphi_{0 \, *} \Cal O_{U_0} =0.
$$
Then there exists a proper surjective morphism
$\varphi : U \to V$ to some normal 4-fold $V$
which factors through $f$ such that

\flushpar
(1) \, \, \, The central fiber of the structure morphism
$V \to \Delta$ is isomorphic to $V_0$, and

\flushpar
(2) \, \, \, $\varphi|_{U_0} = \varphi_0$ under this
identification.
\endproclaim

\proclaim{4.3}\ \ Proof of Theorem 4.1. (Following
the idea of J.Koll\'ar)
\endproclaim

First by the condition (C),
$\biggl( \sum\limits_{k \geq 2} D_k \, . \, R \biggr) =0$,
{\it i.e.\/}

\flushpar
(4.3.1) \qquad \qquad \qquad \qquad \qquad \qquad \, \,
$(D_1 \, . \, R)=0$.

\flushpar
Hence for any irreducible curve $C \subset E$,
we have

\flushpar
(4.3.2) \qquad \qquad \qquad \qquad \quad \, \, \, \,
$(-K_{D_1} \, . \, C)= (-K_X \, . \, C)$.

Let $U$ be an analytic neighborhood of
$E$ in $X$ so that

\flushpar
(4.3.3) \, \, \, $U \cap f^{-1}(t)$ contains
no proper 1-dimensional subspaces \,
$(\forall t \in \Delta - \{0\})$,

\flushpar
and consider $g|_U : U \to g(U)$.

\flushpar
(4.3.4) \, \, \,
Let $D_1^{\circ} := D_1 \cap U (= D \cap U)$, and let
$g(D_1^{\circ})\sptilde$ be the normalization
of $g(D_1^{\circ})$.

Then the induced morphism
$h^{\circ} : D_1^{\circ} \to g(D_1^{\circ})\sptilde$
is a bimeromorphic morphism such that

\flushpar
(4.3.5) \, \, \,
$-K_{D_1^{\circ}}$ is $h^{\circ}\text{-ample}$.

\flushpar
Consider $\overline{NE} (D_1^{\circ}/g(D_1^{\circ})\sptilde)$.
Choose any extremal ray of \, it, and
let $\varphi_0 : D_1^{\circ} \to \varphi_0(D_1^{\circ})$
be the associated contraction.
Since $\text{Exc } \varphi_0 \subset E$ and
$\dim \text{Exc } \varphi_0 = 2$ (Remark 0.6),

\flushpar
(4.3.6) \qquad \qquad \qquad \qquad \quad \, \, \, \, \,
$\text{Exc } \varphi_0 = E_i$ \, for some $i$.

\flushpar
Let
$$
Q := \varphi_0 (E_i).
$$
Then by $R^1 \varphi_{0 \, *} \Cal O_{D_1^{\circ}} =0$
and Theorem 4.2 we have a proper surjective morphism
$\varphi : U \to V$ to some normal 4-fold $V$ over $\Delta$
such that
$$
\cases
(V)_0 \simeq \varphi_0(D_1^{\circ}), \, \text{ and} \\
\varphi|_{D_1^{\circ}} = \varphi_0.
\endcases
\tag 4.3.7
$$
By (4.3.2), (4.3.3) and (4.3.6),

\flushpar
(4.3.8) \, \, \,
$\varphi$ is a flipping contraction with \,
$\text{Exc } \varphi = E_i$.

\flushpar
By the above construction,

\flushpar
(4.3.9) \, \, \, $((V)_0, Q)$ is a germ of a 3-dimensional
terminal singularity which is a Cartier divisor of the
4-dimensional flipping singularity $(V,Q)$.
Note that $V-Q$ is Gorenstein, since so is $U$.

If the singularity index of $((V)_0, Q)$ is greater than 1,
then this is never deformed to be Gorenstien
[Sc], [R1,4], [Mo3,4], [KS],
a contradiction to (4.3.9). So
$((V)_0, Q)$ must be Gorenstein, and hence
a hypersurface singularity [loc.cit].
Then $(V,Q)$ is also a hypersurface singularity (4.3.9),
which again contradicts (4.3.9), since flipping singularities
can never be $\Bbb Q\text{-Gorenstein}$.
Hence the Theorem 4.1.
\quad \qed

By completely the same argument as in 4.3,
we can prove the following which might be a little stronger
than Theorem 4.1:

\proclaim{Remark 4.4}\ \
Let $g : X \to Y$ be as in Definition 0.3, and
let $E$ be any connected component of $\text{Exc } g$ as in
Main Theorem 0.5. Instead of (0.5.1), we assume the
following:

\flushpar
(4.4.1) \, \, \, $E \subset D_1$, and
$g(D_1)$ has a terminal singularity at $g(E)$.

\flushpar
Then

\flushpar
\qquad \qquad \qquad \qquad \qquad \qquad \qquad \, \, \,
$E \cap \bigcup\limits_{k \geq 2} D_k \not= \emptyset$.
\quad \qed
\endproclaim

\head \S 5.\ Description of flips.
\endhead

\definition{Notation 5.0}\ \
Let $g : X \to Y$ be as in Main Theorem 0.5.
$g$ is of type either (A-I), (A-II) or (B).
Let $g^+ : X^+ \to Y$ be the flip of $g$,
let $E^+ \simeq \Bbb P^1$ be the flipped curve,
and $D_k^+$ the proper transform of $D_k$
$(k = 1, \dots , r)$. Moreover, let
$$
\eta := g^{+ \, -1} \circ g : X \dashrightarrow X^+.
$$
Recall that
$N_{E^+/X^+} \simeq \Cal O_{\Bbb P^1}(-1)^{\oplus 3}$ [Kaw4].
\enddefinition

In the following 5.1 through 5.3,
we shall describe the flip $g^+ : X^+ \to Y$,
depending on the type of $g$.

\proclaim{5.1}\ \ (For the Type (A-I))

If $g : X \to Y$ is of Type (A-I),
then

\flushpar
(1) \qquad \qquad \quad
$E^+ \cap D_1^+ \cap D_2^+ = \emptyset$,
and
$(D_1^+ \, . \, E^+) = (D_2^+ \, . \, E^+) = 1$.

\flushpar
(2) \, \, \,
$E^+ \subset D_3^+ \cap D_4^+$,
$D_3^+$, $D_4^+$ are smooth along $E^+$, and
$$
N_{E^+ / D_3^+} \simeq N_{E^+ / D_4^+} \simeq
\Cal O_{\Bbb P^1}(-1)^{\oplus 2}.
$$

\flushpar
(3) \qquad \qquad \qquad \qquad \,
$N_{D_3 \cap E / D_3} \simeq N_{D_4 \cap E / D_4} \simeq
\Cal O_{\Bbb P^1}(-1)^{\oplus 2}$,

\flushpar
and $\eta |_{D_k} : D_k \dashrightarrow D_k^+$
gives the flop of $D_k \supset D_k \cap E$
$(k = 3,4)$.
\endproclaim

\proclaim{5.2}\ \ (For the Type (A-II))

If $g : X \to Y$ is of Type (A-II),
then

\flushpar
(1) \qquad \qquad \quad
$E^+ \cap D_1^+ \cap D_2^+ = \emptyset$,
and
$(D_1^+ \, . \, E^+) = (D_2^+ \, . \, E^+) = 1$.

\flushpar
(2) \, \, \,
$E^+ \subset D_3^+$, and
$D_3^+$ has the canonical singularity of type $\dfrac{1}{2} (1,1,0)$
along $E^+$.

\flushpar
(3) \qquad \qquad \qquad \qquad \qquad \quad
$N_{D_3 \cap E / D_3} \simeq \Cal O_{\Bbb P^1}(-2)^{\oplus 2}$,

\flushpar
and $\eta |_{D_3} : D_3 \dashrightarrow D_3^+$
gives the anti-flip of $D_3 \supset D_3 \cap E$.
\endproclaim

\proclaim{5.3}\ \ (For the Type (B))

If $g : X \to Y$ is of Type (B),
then

\flushpar
(1) \qquad \qquad \qquad \qquad \qquad \qquad \quad \,
$(D_1^+ \, . \, E^+) = 1$.

\flushpar
(2) \, \, \,
$E^+ \subset D_2^+$, and
$D_2^+$ is smooth along $E^+$.

\flushpar
(3) \qquad \qquad \qquad \qquad \qquad \quad
$N_{D_2 \cap E / D_2} \simeq \Cal O_{\Bbb P^1}(-1)^{\oplus 2}$,

\flushpar
and $\eta |_{D_2} : D_2 \dashrightarrow D_2^+$
gives the flop of $D_2 \supset D_2 \cap E$.
\endproclaim

\definition{Remark 5.4}\ \
If $g$ is of Type (A-I) or (B), then
the flip $g^+$ again satisfies all
the conditions of the Assumption 0.2,

On the other hand, if $g$ is of Type (A-II), then
$X^+$ is still smooth while $D_3^+$ is singular,
so (0.2.1) fails for such $g^+$.

This is the reason why in (0.3.1) we eliminated the condition
(0.2.1) from the definition of flips.
\enddefinition

\enddocument